\documentclass[amstex,12pt,russian,amssymb]{article}

\usepackage{mathtext}
\usepackage[cp1251]{inputenc}
\usepackage[T2A]{fontenc}
\usepackage[russian]{babel}
\usepackage[dvips]{graphicx}
\usepackage{amsmath}
\usepackage{amssymb}
\usepackage{amsxtra}
\usepackage{latexsym}
\usepackage{ifthen}

\begin{document}

\newcounter{lemma}
\newcommand{\lemma}{\par \refstepcounter{lemma}%
{\bf Лемма \arabic{lemma}.}}

\newcounter{corollary}
\newcommand{\corollary}{\par \refstepcounter{corollary}%
{\bf Следствие \arabic{corollary}.}}

\newcounter{remark}
\newcommand{\remark}{\par \refstepcounter{remark}%
{\bf Замечание \arabic{remark}.}}

\newcounter{theorem}
\newcommand{\theorem}{\par \refstepcounter{theorem}%
{\bf Теорема \arabic{theorem}.}}

\newcounter{proposition}
\newcommand{\proposition}{\par \refstepcounter{proposition}%
{\bf Предложение \arabic{proposition}.}}

\renewcommand{\refname}{\centerline{\bf Список литературы}}

\newcommand{\proof}{{\it Доказательство.\,\,}}

\noindent УДК 517.5

{\bf Е.А.~Севостьянов} (Житомирский государственный университет
имени Ивана Франко)

{\bf Є.О.~Севостьянов} (Житомирський державний університет імені
Івана Фран\-ка)

{\bf E.A.~Sevost'yanov} (Zhytomyr Ivan Franko State University)

\medskip
{\bf О  существовании решения уравнения Бельтрами с вырождением}

{\bf Про існування розв'язку рівняння Бельтрамі з виродженням }

{\bf On the existence of a solution of the Beltrami equation with
degeneration}

\medskip\medskip
Найдено одно из возможных условий, при которых уравнение Бельтрами с
вырождением эллиптичности имеет непрерывное решение класса Соболева.
При некоторых дополнительных требованиях указанное решение является
гомеоморфным.

\medskip\medskip
Знайдено одну з можливих умов, при якій рівняння Бельтрамі з
виродженням еліптичності має неперервний розв'язок класу Соболєва.
За певних додаткових вимог вказаний розв'язок є гомеоморфним.

\medskip\medskip
We have found one of the possible conditions under which the
Beltrami equation with degeneration of ellipticity has a continuous
solution of the Sobolev class. With some additional requirements,
this solution is homeomorphic.

\newpage
{\bf 1. Введение.} В последние годы активно развивалась тематика,
связанная с существованием решений вырожденных дифференциальных
уравнений Бельтрами, см., напр., \cite{RSY$_1$}, \cite{RSY$_2$},
\cite{GRY} и~\cite{GRSY}. Основные результаты на эту тему собраны в
относительно свежей монографии~\cite{GRSY}, где имеются ссылки на
публикации этих и других авторов. Одна из задач, стоящих при
исследовании уравнений Бельтрами, состоит в поиске условий на
комплексный коэффициент, обеспечивающих существование решений этих
уравнений. Поиск решений, как правило, осуществляется в классе
$ACL$-гомеоморфизмов, хотя вполне корректно рассматривать в этом
качестве и просто непрерывные $ACL$-решения. В данной заметке
получен ещё один результат о существовании решений уравнения
Бельтрами с вырождением, который основан на переходе к обратным
отображениям. В сравнении с работами~\cite{RSY$_1$}, \cite{RSY$_2$}
и \cite{GRY}, мы ослабляем условия на комплексный коэффициент,
требуя только его интегрируемость и не прибегая к ограничениям более
специального вида. Полученное решение уравнения может оказаться не
гомеоморфным, однако, по отношению к предшествующим результатам
степень его гладкости несколько выше и соответствует классу~$W_{\rm
loc}^{1, 2}.$

\medskip
Обратимся к определениям. Всюду далее отображение $f:D\rightarrow
{\Bbb C}$ области $D\subset{\Bbb C}$ предполагается {\it сохраняющим
ориентацию,} в частности, если $f$ -- гомеоморфизм и $z\in D$ --
какая либо его точка дифференцируемости, то {\it якобиан} этого
отображения в точке $z$ положителен. Для комплекснозначной функции
$f:D\rightarrow {\Bbb C},$ заданной в области $D\subset {\Bbb C},$
имеющей частные производные по $x$ и $y$ при почти всех $z=x + iy,$
полагаем $f_{\overline{z}} = \left(f_x + if_y\right)/2$ и $f_z =
\left(f_x - if_y\right)/2.$ {\it Комплексной дилатацией} отображения
$f$ в точке $z$ называется функция $\mu:D\rightarrow {\Bbb C},$
определённая равенством $\mu(z)=\mu_f(z)=f_{\overline{z}}/f_z,$ при
$f_z \ne 0$ и $\mu(z)=0$ в противном случае. {\it Мак\-си\-маль\-ной
дилатацией} отображения $f$ в точке $z$ называется следующая
функция:
\begin{equation}\label{eq1}
K_{\mu}(z)=K_{\mu_f}(z)=\quad\frac{1+|\mu (z)|}{1-|\mu\,(z)|}\,.
\end{equation}
Если задана измеримая по Лебегу функция $\mu:D\rightarrow {\Bbb D},$
${\Bbb D}=\{z\in {\Bbb C}: |z|<1\},$ то без привязки к какому-либо
отображению $f$ будем называть величину, вычисляемую при помощи
равенства~(\ref{eq1}), максимальной дилатацией, соответствующей
функции $\mu.$ Заметим, что якобиан отображения $f$ в точке $z\in D$
может быть вычислен при помощи равенства
$$J(z,
f)=|f_z|^2-|f_{\overline{z}}|^2\,,$$
что может быть проверено прямым подсчётом. Нетрудно видеть, что
$K_{\mu_f}(z)=\frac{|f_z|+|f_{\overline{z}}|}{|f_z|-|f_{\overline{z}}|}$
во всех точках $z\in D$ отображения $f,$ имеющего частные
производные в точке $z,$ где якобиан $J(z, f)$ не обращается в нуль.
Напомним, что отображение $f:D\rightarrow {\Bbb C}$ называется {\it
квазиконформным}, если $f$ -- гомеоморфизм класса~$W_{\rm loc}^{1,
2}(D)$ и, кроме того, найдётся постоянная $K\geqslant 1$ такая, что
$\Vert f^{\,\prime}(z)\Vert^2\leqslant K\cdot |J(z, f)|,$ где $\Vert
f^{\,\prime}(z)\Vert=|f_z|+|f_{\overline{z}}|.$ {\it Уравнением
Бельтрами} будем называть функциональное уравнение вида
\begin{equation}\label{eq2}
f_{\overline{z}}=\mu(z)\cdot f_z\,,
\end{equation}
в котором $\mu=\mu(z)$ -- заданная неизвестная функция. Для
фиксированного натурального числа $k\geqslant 1$ обозначим
\begin{equation}\label{eq12:} \mu_k(z)= \left
\{\begin{array}{rr}
 \mu(z),& K_{\mu}(z)\leqslant k,
\\ 0\ , & K_{\mu}(z)> k\,.
\end{array} \right.
\end{equation}
Пусть $f_k$ -- гомеоморфное $ACL$-решение уравнения
$f_{\overline{z}}=\mu_k(z)\cdot f_z,$ отображающее единичный круг на
себя, удовлетворяющее условиям нормировки $f(0)=0,$ $f(1)=1$
существующее ввиду~\cite[теорема~8.2]{Bo}. Пусть $g_k$ -- обратное
отображение к $f_k,$ тогда его комплексная дилатация $\mu_{g_k}$
вычисляется согласно соотношению
$\mu_{g_k}(w)=-\mu_k(g_k(w))=-\mu_k(f^{\,-1}_k(w)),$ см.
напр.,~\cite[(4).C.I]{A}. Тогда максимальная дилатация отображения
$g_k$ вычисляется по соотношению
\begin{equation}\label{eq3}
K_{\mu_{g_k}}(w)=\quad\frac{1+|\mu_k(f^{\,-1}_k(w))|}{1-|\mu_k(f^{\,-1}_k(w))|}\,.
\end{equation}

\medskip
Будем говорить, что функция ${\varphi}:D\rightarrow{\Bbb R}$ имеет
{\it конечное среднее колебание} в точке $x_0\in D$, пишем
$\varphi\in FMO(x_0),$ если
\begin{equation}\label{eq17:}
{\limsup\limits_{\varepsilon\rightarrow
0}}\frac{1}{\Omega_n\varepsilon^n}\int\limits_{B(
x_0,\,\varepsilon)}
|{\varphi}(x)-\overline{{\varphi}}_{\varepsilon}|\ dm(x)\, <\,
\infty\,,
\end{equation}
где $\Omega_n$ -- объём единичного шара в ${\Bbb R}^n,$
$\overline{{\varphi}}_{\varepsilon}=\frac{1}{\Omega_n\varepsilon^n}\int\limits_{B(
x_0,\,\varepsilon)} {\varphi}(x)\ dm(x),$ см., напр.,
\cite[разд.~2]{RSY$_2$}.
Заметим, что при выполнении условия $(\ref{eq17:})$ возможна
ситуация, когда
$\overline{{\varphi}_{\varepsilon}}\rightarrow\infty$ при
$\varepsilon\rightarrow 0.$
Также будем говорить, что  ${\varphi}:D\rightarrow{\Bbb R}$ --
функция конечного среднего колебания \ в \ области\  D, пишем
${\varphi}\in FMO(D),$ если ${\varphi}$ имеет конечное среднее
колебание в каждой точке  $x_0\in D.$ Справедливо следующее
утверждение.

\medskip
\begin{theorem}\label{th1}{\sl\,
Пусть $Q:{\Bbb D}\rightarrow[1, \infty]$ -- интегрируемая в ${\Bbb
D}$ функция и пусть функция $\mu:{\Bbb D}\rightarrow {\Bbb D}$
измерима по Лебегу. Предположим, что для почти всех $w\in {\Bbb D}$
\begin{equation}\label{eq10}
K_{\mu_{g_k}}(w)\leqslant Q(w)\,,
\end{equation}
где $g_k=f_k^{\,-1}$ и $f_k$ -- гомеоморфное $ACL$-решение уравнения
$f_{\overline{z}}=\mu_k(z)\cdot f_z,$ отображающее единичный круг на
себя, удовлетворяющее условиям нормировки $f(0)=0,$ $f(1)=1,$ кроме
того, $\mu_k(z)$ задаётся соотношением~(\ref{eq12:}).

Тогда уравнение~(\ref{eq2}) имеет непрерывное $W_{\rm loc}^{1,
2}({\Bbb D})$-решение $f:{\Bbb D}\rightarrow {\Bbb D},$
удовлетворяющее условию
\begin{equation}\label{eq11}
|f(z)-f(z_0)|\leqslant\frac{C\cdot (\Vert
Q\Vert_1)^{1/2}}{\log^{1/2}\left(1+\frac{r_0}{|z-z_0|}\right)}\quad\forall\,\,z\in
B(z_0, r_0)
\end{equation}
в произвольной точке $z_0\in {\Bbb D},$ где $\Vert Q\Vert_1$ --
норма $Q$ в $L^1({\Bbb D}),$ $C$ -- некоторая постоянная и
$0<2r_0<{\rm dist}\,(z_0,
\partial {\Bbb D})$ -- произвольно.
Если дополнительно $Q(z)\in FMO({\Bbb D}),$ либо
\begin{equation}\label{eq5**}
\int\limits_{0}^{\delta(w_0)}\frac{dt}{tq_{w_0}(t)}=\infty
\end{equation}
для каждого $w_0\in {\Bbb D}$ и некотором $\delta(w_0)>0,$
$q_{w_0}(r)=\frac{1}{2\pi}\int\limits_{0}^{2\pi}Q(w_0+re^{\,i\theta})\,d\theta,$
то $f$ является гомеоморфизмом в ${\Bbb D}.$ }
\end{theorem}

\medskip
{\bf 2. Основная лемма о сходимости.} Связь между сходимостью
отображений и поведением их комплексных коэффициентов является
важнейшим элементом, используемым при доказательстве основной
теоремы. По поводу аналогичных утверждений, известных на данный
момент, укажем, напр., на~\cite[гл.~2]{GRSY} либо~\cite{RSY}. Что
касается изучаемого в работе случая, справедлива следующая лемма.

\medskip
\begin{lemma}\label{lem1}
{\sl\, Пусть $Q\in L^1(D),$ $\mu:D\rightarrow {\Bbb D}$ -- измеримая
по Лебегу функция, и пусть $f_k,$ $k=1,2,\ldots $ --
последовательность сохраняющих ориентацию гомеоморфизмов области $D$
на себя, принадлежащих классу $W_{\rm loc}^{1, 2}(D)$ и имеющих
комплексные коэффициенты $\mu_k(z).$ Предположим, что $f_k$ сходится
локально равномерно в $D$ к отображению $f:D\rightarrow {\Bbb C},$ а
последовательность $\mu_k(z)$ сходится к $\mu$ при $k\rightarrow
\infty$ при почти всех $z\in D.$ Пусть также обратные отображения
$g_k:=f_k^{\,-1}$ принадлежат классу $W_{\rm loc}^{1, 2}(D),$ при
этом, при почти всех $w\in D$
$$K_{\mu_{g_k}}(w)\quad=\quad\frac{1+|\mu_k(f^{\,-1}_k(w))|}{1-|\mu_k(f^{\,-1}_k(w))|}\leqslant Q(w)\,.$$
Тогда $f\in W_{\rm loc}^{1, 2}(D)$ и $\mu$ -- комплексная
характеристика отображения $f,$ т.е., $f_{\overline{z}}=\mu(z)\cdot
f_z$ при почти всех $z\in D.$
 }
\end{lemma}

\medskip
\begin{proof}
Будем в целом следовать схеме, изложенной при
доказательстве~\cite[теорема~3.1]{RSY}, см.
также~\cite[теорема~2.1]{GRSY} и \cite[лемма~III.3.5]{Re}. Обозначим
$\partial f=f_z$ и $\overline{\partial}f=f_{\overline{z}}.$ Пусть
$C$ -- произвольный компакт в $D.$ Поскольку по предположению
отображения~$g_k=f_k^{\,-1}$ принадлежат классу $W_{\rm loc}^{1,
2},$ то $g_k$ обладают $N$-свойством Лузина, см., напр.,
\cite[следствие~B]{MM}. Тогда якобиан $J(z, f)$ почти всюду отличен
от нуля, см., напр., \cite[теорема~1]{Pon}, более того, имеет место
замена переменных в интеграле, см., напр.,~\cite[теорема~3.2.5]{Fe}.
В таком, случае, будем иметь:
$$\int\limits_{C}|\partial f_k(z)|^2\,dm(z)=
\int\limits_C\left(|\partial
f_k(z)|^{\,2}-|\overline{\partial}f_k(z)|^{\,2}\right)\cdot
\frac{|\partial f_k(z)|^{\,2} dm(z)}{\left(|\partial
f_k(z)|^{\,2}-|\overline{\partial}f_k(z)|^{\,2}\right)}=$$
\begin{equation}\label{eq4}
=\quad \int\limits_C J(z,
f_k)\cdot\frac{1}{1-\left|\frac{\overline{\partial}f_k(z)}{\partial
f_k(z)}\right|^2}\,\,dm(z)= \quad\int\limits_{f_k(C)}
\frac{dm(w)}{1-|\mu_k(f_k^{\,-1}(w))|^2}\quad\leqslant
\end{equation}
$$\leqslant \quad\int\limits_{D} Q(w)\, dm(w)<\infty\,.$$
Из~(\ref{eq4}) вытекает, что $f\in W_{\rm loc}^{1, 2},$ при этом,
$\partial f_k$ и $\overline{\partial} f_k$ слабо сходятся в $L_{\rm
loc}^1(D)$ к $\partial f$ и $\overline{\partial} f,$ соответственно
(см.~\cite[теорема~3.1]{RSY} и \cite[лемма~III.3.5]{Re}).

\medskip
Осталось показать, что отображение $f$ является решением уравнения
Бельтрами $f_{\overline{z}}=\mu(z)\cdot f_z.$ Положим
$\zeta(z)=\overline{\partial} f(z)-\mu(z)\partial f(z)$ и покажем,
что $\zeta(z)=0$ почти всюду. Пусть $B$ -- произвольный круг,
лежащий вместе со своим замыканием в ${\Bbb D}.$ По неравенству
треугольника
\begin{equation}\label{eq9}
\left|\int\limits_B\zeta(z)\,dm(z)\right|\leqslant I_1(k)+I_2(k)\,,
\end{equation}
где
\begin{equation}\label{eq7}
I_1(k)=\left|\int\limits_B(\overline{\partial}f(z)-\overline{\partial}f_k(z))\,dm(z)\right|
\end{equation}
и
\begin{equation}\label{eq8}
I_2(k)=\left|\int\limits_B(\mu(z)\partial f(z)-\mu_k(z)\partial
f_k(z))\,dm(z)\right|\,.
\end{equation}
Ввиду доказанного выше, $I_1(k)\rightarrow 0$ при
$k\rightarrow\infty.$ Осталось разобраться с выражением $I_2(k).$
Для этого заметим, что
$I_2(k)\leqslant I^{\,\prime}_2(k)+I^{\,\prime\prime}_2(k),$ где
$$I^{\,\prime}_2(k)=
\left|\int\limits_B\mu(z)(\partial f(z)-\partial
f_k(z))\,dm(z)\right|$$
и
$$I^{\,\prime\prime}_2(k)=
\left|\int\limits_B(\mu(z)-\mu_k(z))\partial f_k(z)\,dm(z)\right|$$
Ввиду слабой сходимости $\partial f_k\rightarrow \partial f$ в
$L^1_{\rm loc}(D)$ при $k\rightarrow\infty,$ мы получим, что
$I^{\,\prime}_2(k)\rightarrow 0$ при $k\rightarrow\infty,$ поскольку
$\mu\in L^{\infty}(D).$ Более того, для заданного $\varepsilon>0$
отыщется $\delta=\delta(\varepsilon)>0$ такое, что, как только
$m(E)<\delta,$ $E\subset B,$ то
\begin{equation}\label{eq5}\int\limits_E|\partial f_k(z)|\,dm(z)\leqslant
\int\limits_E|\partial f_k(z)-\partial f(z)|\,dm(z)+
\int\limits_E|\partial f(z)|\,dm(z)<\varepsilon\,,
\end{equation}
$k=1,2\ldots ,$ поскольку $\partial f_k\rightarrow \partial f$ слабо
в $L^1_{\rm loc}({\Bbb D}),$ кроме того, отображение $\partial f$
интегрируемо с квадратом по доказанному выше, а значит, имеет место
абсолютная непрерывность интеграла Лебега.

\medskip
Окончательно, по теореме Егорова (см.~\cite[теорема~III.6.12]{Sa})
для каждого $\delta>0$ найдётся множество $S\subset B$ такое, что
$m(B\setminus S)<\delta$ и $\mu_k(z)\rightarrow \mu(z)$ равномерно
на~$S.$ Тогда $|\mu_k(z)-\mu(z)|<\varepsilon$ при всех $k\geqslant
k_0,$ некотором $k_0=k_0(\varepsilon)$ и всех $z\in S,$ кроме того,
ввиду~(\ref{eq5}), а также ввиду~(\ref{eq4}) и по неравенству
Гёльдера мы имеем, что
$$I^{\,\prime\prime}_2(k)\leqslant \varepsilon \int\limits_S|\partial
f_k(z)|\,dm(z)+ 2\int\limits_{B\setminus S}|\partial
f_k(z)|\,dm(z)<$$
\begin{equation}\label{eq6}
<\varepsilon\cdot\left\{\left(\int\limits_D Q(w)\,
dm(w)\right)^{1/2}\cdot (m(C))^{1/2}+2\right\}
\end{equation}
при тех же $k\geqslant k_0.$ Из~(\ref{eq9}), (\ref{eq7}),
(\ref{eq8}) и (\ref{eq6}) вытекает, что
$\int\limits_B\zeta(z)\,dm(z)=0$ для всех кругов $B,$ компактно
вложенных в ${\Bbb D}.$ на основании теоремы Лебега о
дифференцировании неопределённого интеграла (см.
\cite[IV(6.3)]{Sa}), отсюда вытекает, что $\zeta(z)=0$ почти всюду в
$D.$ Лемма доказана.~$\Box$
\end{proof}

\medskip
{\bf 3. Доказательство теоремы~\ref{th1}}. Пусть вначале функция
$Q,$ заданная по условию теоремы, просто интегрируема в ${\Bbb D}.$
Рассмотрим по\-сле\-до\-ва\-тель\-ность комплекснозначных функций
\begin{equation}\label{eq12}
\mu_k(z)= \left \{\begin{array}{rr}
 \mu(z),& K_{\mu}(z)\leqslant k,
\\ 0\ , & K_{\mu}(z)>k,
\end{array} \right.
\end{equation}
где $K_{\mu}(z)$ определяется соотношением~(\ref{eq1}). Заметим, что
$\mu_k(z)\leqslant\frac{k-1}{k+1}<1,$ поэтому уравнение
$(\ref{eq2}),$ где вместо $\mu$ в правой части взято $\mu:=\mu_k,$ а
$\mu_k$ определено соотношениями~(\ref{eq12}), имеет гомеоморфное
$W_{\rm loc}^{1, 2}({\Bbb D})$-решение $f_k:{\Bbb D}\rightarrow{\Bbb
C}$ с нормировками $f_k(0)=0,$ $f_k(1)=1,$ которое является
$k$-квазиконформным в ${\Bbb D},$ см.~\cite[теорема~8.2]{Bo}. Ввиду
этой же теоремы, $f_k$ отображают единичный круг на себя, при этом,
$g_k=f^{\,-1}_k$ также являются квазиконформными, в частности,
принадлежат класу $W_{\rm loc}^{1, 2}({\Bbb D})$
(см.~\cite[теорема~9.1]{BI}). В силу~\cite[теорема~6.10]{MRSY$_1$} и
ввиду условия~(\ref{eq10}) для каждого $k\in {\Bbb N}$
\begin{equation} \label{eq2*B}
M(g_k(\Gamma))\leqslant \int\limits_{\Bbb
D}K_{\mu_{g_k}}(w)\cdot\rho_*^2 (w) \,dm(w)\leqslant
\int\limits_{\Bbb D}Q(w)\cdot\rho_*^2 (w) \,dm(w)
\end{equation}
для произвольного семейства кривых $\Gamma$ в ${\Bbb D}$ и каждой
функции $\rho_*\in {\rm adm}\,\Gamma,$ где $M$ -- модуль семейства
кривых (см., напр.,~\cite[разд.~6]{Va}). В
силу~\cite[теорема~1.1]{SevSkv} семейство отображений $f_k$
равностепенно непрерывно в ${\Bbb D}.$ Значит, ввиду теоремы
Арцела-Асколи $f_k$ является нормальным семейством
(см.~\cite[теорема~20.4]{Va}), другими словами, найдётся
подпоследовательность $f_{k_l}$ последовательности $f_k,$ сходящаяся
локально равномерно в ${\Bbb D}$ к некоторому отображению $f:{\Bbb
D}\rightarrow \overline{{\Bbb D}}.$ Заметим также, что
$\mu_k(z)\rightarrow \mu(z)$ при $k\rightarrow\infty$ для почти всех
$z\in {\Bbb D},$ поскольку $|\mu(z)|<1$ и, значит, $K_{\mu}(z)$
в~(\ref{eq1}) конечна при всех $z\in {\Bbb D}.$ Тогда по
лемме~\ref{lem1} отображение $f$ принадлежит классу $W_{\rm loc}^{1,
2}({\Bbb D})$ и является решением исходного уравнения
Бельтрами~(\ref{eq2}).

\medskip
По~\cite[теорема~1]{SSD}
$$|f_k(z)-f_k(z_0)|\leqslant\frac{C\cdot (\Vert
Q\Vert_1)^{1/2}}{\log^{1/2}\left(1+\frac{r_0}{|z-z_0|}\right)}\quad\forall\,\,z\in
B(z_0, r_0)$$
в произвольной точке $z_0\in {\Bbb D},$ где $\Vert Q\Vert_1$ --
норма $Q$ в $L^1({\Bbb D}),$ $C$ -- некоторая постоянная и
$0<2r_0<{\rm dist}\,(z_0,
\partial {\Bbb D}).$ Переходя здесь к пределу при
$k\rightarrow\infty,$ получаем соотношение~(\ref{eq11}). Первая
часть утверждения теоремы~\ref{th1} установлена.

\medskip
Предположим теперь, что~$Q\in FMO({\Bbb D}),$ либо выполняется
соотношение~(\ref{eq5**}). Тогда последовательность $g_k$ образует
равностепенно непрерывное семейство отображений
(см.~\cite[теоремы~6.1 и 6.5]{RS}). Значит, ввиду теоремы
Арцела-Асколи $g_k$ является нормальным семейством
(см.~\cite[теорема~20.4]{Va}), другими словами, найдётся
подпоследовательность $g_{k_l}$ последовательности $g_k,$ сходящаяся
локально равномерно в ${\Bbb D}$ к некоторому отображению $g:{\Bbb
D}\rightarrow \overline{{\Bbb D}}.$ В силу условия нормировки
$g_{k_l}(0)=0$ и $g_{k_l}(1)=1$ при всех $l=1,2,\ldots .$ Тогда в
силу~\cite[теорема~4.1]{RSS} отображение $g$ является гомеоморфизмом
в ${\Bbb D},$ кроме того, по~\cite[лемма~3.1]{RSS} мы имеем также,
что $f_{k_l}\rightarrow f=g^{\,-1}$ при $l\rightarrow\infty$
локально равномерно в ${\Bbb D}.$ Далее применяем схему рассуждений,
уже применённую выше к случаю интегрируемой функции $Q.$ Поскольку
$\mu_k(z)\rightarrow \mu(z)$ при $k\rightarrow\infty$ и при почти
всех $z\in {\Bbb D},$ то по лемме~\ref{lem1} отображение $f$
принадлежит классу $W_{\rm loc}^{1, 2}({\Bbb D})$ и является
решением исходного уравнения Бельтрами~(\ref{eq2}). Теорема
доказана.~ $\Box$

\medskip
{\bf Пример.} Пусть $p\geqslant 1$ -- произвольное число и пусть
$0<\alpha<2/p.$ Как обычно, мы используем запись $z=re^{i\theta},$
$r\geqslant 0$ и $\theta\in [0, 2\pi).$ Положим
\begin{equation}\label{eq3A}\mu(z)= \left
\{\begin{array}{rr}
 e^{2i\theta}\frac{2r-\alpha(2r-1)}{2r+\alpha(2r-1)},& 1/2<|z|<1,
\\ 0\ , & |z|\leqslant 1/2\,.
\end{array} \right.
\end{equation}
Используя соотношение
$$\mu_f(z)=\frac{\overline{\partial} f}{\partial f}=e^{2i\theta}\frac{rf_r+if_{\theta}}{rf_r-if_{\theta}}\,,$$
см. равенство~(11.129) в \cite{MRSY}, мы получаем, что отображение
\begin{equation}\label{eq4A}f(z)=\left
\{\begin{array}{rr}
 \frac{z}{|z|}(2|z|-1)^{1/\alpha},& 1/2<|z|<1,
\\ 0\ , & |z|\leqslant 1/2
\end{array} \right.
\end{equation}
является решением уравнения~$f_{\overline{z}}=\mu(z)\cdot f_z,$ где
функция $\mu$ задаётся соотношением~(\ref{eq3A}). Заметим, что
существование решения данного уравнения обеспечивается
теоремой~\ref{th1} (для этого проверим выполнение всех условий этой
теоремы). Заметим, что для заданной соотношением~(\ref{eq3A})
функции $\mu$ соответствующей ей максимальной дилатацией $K_{\mu}$
будет функция
\begin{equation}\label{eq5A}K_{\mu}(z)=\left
\{\begin{array}{rr}
 \frac{4|z|}{2\alpha(2|z|-1)},& 1/2<|z|<1,
\\ 1\ , & |z|\leqslant 1/2
\end{array} \right.\,.
\end{equation}
Заметим, что $K_{\mu}(z)\leqslant k$ при $|z|\geqslant
\frac{1}{2}\cdot \frac{k\alpha}{k\alpha-1}$ и $K_{\mu}(z)>k$ в
противном случае. Пусть, как и прежде,
\begin{equation*}\label{eq5B}\mu_k(z)= \left
\{\begin{array}{rr}
 \mu(z),& K_{\mu}(z)\leqslant k,
\\ 0\ , & K_{\mu}(z)> k\,.
\end{array} \right.
\end{equation*}
Заметим, что решениями уравнения $f_{\overline{z}}=\mu_k(z)\cdot
f_z$ являются отображения
\begin{equation*}\label{eq6A}f_k(z)=\left
\{\begin{array}{rr}
 \frac{z}{|z|}(2|z|-1)^{1/\alpha},& \frac{1}{2}\cdot \frac{k\alpha}{k\alpha-1}<|z|<1,
\\ 0\ , & |z|\leqslant \frac{k\alpha}{k\alpha-1}
\end{array} \right.\,,
\end{equation*}
при этом, обратные отображения $g_k(y)=f_k^{\,-1}(y)$ вычисляются по
формуле
\begin{equation}\label{eq7A}g_k(y)=\left
\{\begin{array}{rr}
 \frac{y(|y|^{\alpha}+1)}{2|y|},& |y|>\left(\frac{k\alpha}{k\alpha-1}-1\right)^{1/\alpha},
\\ \frac{y\cdot\frac{k\alpha}{2(k\alpha-1)}}{\left(\frac{k\alpha}{k\alpha-1}-1\right)^{1/\alpha}}
, & |y|\leqslant\left(\frac{k\alpha}{k\alpha-1}-1\right)^{1/\alpha}
\end{array} \right.\,.
\end{equation}
Из~(\ref{eq5A})) вытекает, что
\begin{equation}\label{eq7B}K_{\mu_k}(z)=\left
\{\begin{array}{rr}
 \frac{4|z|}{2\alpha(2|z|-1)},& \frac{1}{2}\cdot \frac{k\alpha}{k\alpha-1}<|z|<1,
\\ 1\ , & |z|\leqslant \frac{k\alpha}{k\alpha-1}
\end{array} \right.\,.
\end{equation}
Нам следует проверить выполнение~(\ref{eq10}) для некоторой
интегрируемой в ${\Bbb D}$ функции $Q.$ Для этой цели, подставим
отображения $g_k$ из~(\ref{eq7A}), в максимальную
дилатацию~$K_{\mu_k},$ определённую равенством~(\ref{eq7B}). Тогда
\begin{equation*}\label{eq8A}K_{\mu_{g_k}}(y)=\left
\{\begin{array}{rr}
 \frac{|y|^{\alpha}+1}{\alpha|y|^{\alpha}},&
 |y|>\frac{y\cdot\frac{k\alpha}{2(k\alpha-1)}}{\left(\frac{k\alpha}{k\alpha-1}-1\right)^{1/\alpha}}\,,
\\ 1\ , & |y|\leqslant\left(\frac{k\alpha}{k\alpha-1}-1\right)^{1/\alpha}
\end{array} \right.\,.
\end{equation*}
Заметим, что $K_{\mu_{g_k}}(y)\leqslant Q(y):=
\frac{|y|^{\alpha}+1}{\alpha|y|^{\alpha}}$ при всех $y\in {\Bbb
B}^n.$ При этом, функция $Q$ интегрируема в ${\Bbb B}^n$ даже в
степени $p,$ а не только в степени $1$ (см. рассуждения,
использованные при рассмотрении~\cite[предложение~6.3]{MRSY}). По
построению $f_k(0)=0$ и $f_k(1)=1.$ Поэтому все условия
теоремы~\ref{th1} выполняются, а в качестве искомого решения
уравнения~$f_{\overline{z}}=\mu_k(z)\cdot f_z$ может быть
рассмотрено отображение $f=f(z),$ определённое
равенством~(\ref{eq4A}). Более того, из доказательства этой теоремы
следует, что отображение $f$ является именно указанным там решением
уравнения, поскольку оно является локально равномерным пределом
последовательности $f_k.$ Отметим, что отображение $f$ не является
гомеоморфным решением, также оно не является открытым и дискретным.

\medskip
Покажем, что для заданной функции $\mu$ гомеоморфного $W_{\rm
loc}^{1, 2}({\Bbb D})$-решения уравнения Бельтрами~(\ref{eq2}) не
существует. В самом деле, пусть $g:{\Bbb D}\rightarrow {\Bbb D}$ --
такое решение. За счёт теоремы Римана об отображении, мы можем
считать, что $g$ отображает единичный круг на себя. Заметим, что при
$1/2<|z|<1$ отображение $f,$ а значит, и отображение $g,$ локально
квазиконформно, поэтому ввиду теоремы единственности $g=\varphi\circ
f,$ где $\varphi$ -- некоторое конформное отображение. Заметим, что
$\varphi$ определено в проколотом круге ${\Bbb D}\setminus\{0\},$
так как $f(\{1/2<|z|<1\})={\Bbb D}\setminus\{0\}.$ Отсюда $g\circ
f^{\,-1}=\varphi,$ и, поскольку $\varphi$ конформно в ${\Bbb
D}\setminus\{0\},$ то оно продолжается по непрерывности в точку $0.$
Но тогда и отображение $g\circ f^{\,-1}$ продолжается по
непрерывности в точку $0,$ что неверно, поскольку $f^{\,-1}(y)=
\frac{y(|y|^{\alpha}+1)}{2|y|},$ а $g$ -- некоторый автоморфизм
единичного круга. Полученное противоречие опровергает предположение
о существовании гомеоморфного решения $g.$


КОНТАКТНАЯ ИНФОРМАЦИЯ

\medskip
\noindent{{\bf Евгений Александрович Севостьянов} \\
{\bf 1.} Житомирский государственный университет им.\ И.~Франко\\
кафедра математического анализа, ул. Большая Бердичевская, 40 \\
г.~Житомир, Украина, 10 008 \\
{\bf 2.} Институт прикладной математики и механики
НАН Украины, \\
отдел теории функций, ул.~Добровольского, 1 \\
г.~Славянск, Украина, 84 100\\
e-mail: esevostyanov2009@gmail.com}

\end{document}